\documentclass[10pt]{amsart}

\usepackage{amssymb}
\usepackage{thmtools}
\usepackage[hidelinks]{hyperref}

\makeatletter
\tagsleft@false
\def\@eqnnum{{\normalfont\normalcolor \tagform@\theequation}}
\let\veqno\@@eqno
\makeatother

\theoremstyle{plain}

\newtheorem{lemma}{Lemma}[section]

\title[Vanishing sets and character degrees]
{Vanishing-set inclusion and character degrees in finite $p$-groups}

\author{Zhi-Lin Zhang}
\address{Independent Researcher, Taipei, Taiwan}
\email{hsa00000@gmail.com}


\begin{document}

\begin{abstract}
Let $G$ be a finite $p$-group and let
$\chi,\psi\in\operatorname{Irr}(G)$. We prove that
\[
  \operatorname{Van}(\chi)\subseteq\operatorname{Van}(\psi)
  \quad\Longrightarrow\quad
  \chi(1)\mid\psi(1),
\]
where
$\operatorname{Van}(\theta)=\{g\in G:\theta(g)=0\}$. Consequently, two
irreducible characters of a finite $p$-group with the same vanishing set
have the same degree. This answers a question attributed to J.~Sangroniz
and recorded by Hung, Moret\'o, and Morotti, extending the previously known
case in which one character has degree $p$. The same divisibility conclusion
holds for irreducible characters of finite nilpotent groups.
\end{abstract}

\maketitle

\section{Introduction}

Throughout the paper, all groups are finite and all characters are ordinary
complex characters. For a character $\theta$ of a group $G$, write
\[
  \operatorname{Van}(\theta)=\{g\in G:\theta(g)=0\}
\]
for its vanishing set. By Burnside's theorem, every nonlinear irreducible
character has a zero \cite[Theorem~3.15]{Isaacs}. Malle, Navarro, and Olsson
proved that such a zero can be chosen on an element of prime-power order
\cite[Theorem~B]{MalleNavarroOlsson}. General background on character zeros
is surveyed in \cite{DolfiPacificiSanus}. For finite $p$-groups, quantitative
questions concerning the number of vanishing conjugacy classes or vanishing
elements were studied in \cite{MoretoSangroniz,MoretoNavarro}. Recent work
on common zeros of two irreducible characters includes \cite{LewisEtAl};
those questions concern nonempty intersections of vanishing sets, whereas
the present paper concerns inclusion of the entire sets.

Hung, Moret\'o, and Morotti isolated a degree-comparison problem for
characters with related vanishing sets. For nonlinear
$\chi,\psi\in\operatorname{Irr}(G)$, the argument of
\cite[Proposition~7.2]{HungMoretoMorotti} shows that coprime degrees
preclude inclusion of the vanishing sets in either direction:
\[
  (\chi(1),\psi(1))=1
  \quad\Longrightarrow\quad
  \operatorname{Van}(\chi)\nsubseteq\operatorname{Van}(\psi)
  \quad\text{and}\quad
  \operatorname{Van}(\psi)\nsubseteq\operatorname{Van}(\chi).
\]
The nonlinearity hypothesis is necessary because every linear character
has empty vanishing set, although that hypothesis is omitted from the
displayed statement of the cited proposition. Thus, for nonlinear
characters, even one-sided inclusion imposes an arithmetic constraint on
the two degrees. They also noted finite groups, including
$\operatorname{SL}_2(5)$ and $\operatorname{PSL}_2(11)$, for which equality
of vanishing sets does not force equality of degrees
\cite[discussion preceding Question~7.3]{HungMoretoMorotti}. By contrast,
Belonogov proved that, for the symmetric group $\mathsf S_n$, the vanishing
set of an irreducible character determines that character up to multiplication
by the sign character \cite[Theorem~1]{Belonogov}.

Hung, Moret\'o, and Morotti recorded Question~7.3, which they attribute to
J.~Sangroniz \cite[Question~7.3]{HungMoretoMorotti}. It asks whether two
irreducible characters of a finite $p$-group with the same vanishing set
must have the same degree. They proved that the answer is affirmative when
one of the characters has degree $p$
\cite[Proposition~7.6]{HungMoretoMorotti}.

A related one-sided comparison was subsequently obtained for symmetric
groups by Giannelli, Law, and McDowell. If
$P\in\operatorname{Syl}_p(\mathsf S_n)$ and
$\chi,\psi\in\operatorname{Irr}(\mathsf S_n)$, write $\chi_P$ and $\psi_P$
for their restrictions to $P$. Then
\[
  \operatorname{Van}(\chi_P)\subseteq\operatorname{Van}(\psi_P)
  \quad\Longrightarrow\quad
  \chi(1)_p\leq\psi(1)_p,
\]
where $m_p$ denotes the $p$-part of a positive integer $m$
\cite[Theorem~S2]{GiannelliLawMcDowell}. They further proved that inclusion
between the zero sets on elements of prime-power order implies the
corresponding inequality between the full character degrees
\cite[Corollary~S3]{GiannelliLawMcDowell}. Their result and ours are
formally parallel but apply in different settings: the restricted characters
$\chi_P$ and $\psi_P$ need not be irreducible, whereas the characters below
are irreducible characters of the $p$-group itself.

Our first headline theorem gives a one-sided answer in the setting of finite
$p$-groups.

\begin{restatable}{theorem}{pgroupdivisibility}
\label{thm:p-group-divisibility}
Let $p$ be a prime, let $G$ be a finite $p$-group, and let
$\chi,\psi\in\operatorname{Irr}(G)$. Then
\[
  \operatorname{Van}(\chi)\subseteq\operatorname{Van}(\psi)
  \quad\Longrightarrow\quad
  \chi(1)\mid\psi(1).
\]
Consequently,
\[
  \operatorname{Van}(\chi)=\operatorname{Van}(\psi)
  \quad\Longrightarrow\quad
  \chi(1)=\psi(1).
\]
\end{restatable}

In particular, Theorem~\ref{thm:p-group-divisibility} answers
\cite[Question~7.3]{HungMoretoMorotti} affirmatively and extends
\cite[Proposition~7.6]{HungMoretoMorotti} from the degree-$p$ case to
arbitrary irreducible-character degrees.

The conclusion extends prime by prime to finite nilpotent groups.

\begin{restatable}{theorem}{nilpotentdivisibility}
\label{thm:nilpotent-divisibility}
Let $G$ be a finite nilpotent group and let
$\chi,\psi\in\operatorname{Irr}(G)$. If
$\operatorname{Van}(\chi)\subseteq\operatorname{Van}(\psi)$, then
$\chi(1)\mid\psi(1)$. In particular, equality of vanishing sets implies
equality of degrees.
\end{restatable}

The main external input is Wilde's restriction-norm calculation
\cite[proof of Lemma~3.3]{Wilde}. Along a fixed chief series of a finite $p$-group,
the restriction norm of an irreducible character changes across each chief
factor by a factor of either $1$ or $p$. The factor is $p$ precisely when
the character vanishes throughout the corresponding set difference in the
series. Using the same chief series for $\chi$ and $\psi$ converts inclusion
of vanishing sets into inclusion of the indices at which these jumps occur.
Theorem~\ref{thm:nilpotent-divisibility} then follows by applying
Theorem~\ref{thm:p-group-divisibility} to the Sylow direct factors.

Section~\ref{sec:p-groups} proves the $p$-group theorem, including the
restriction-norm calculation. Section~\ref{sec:nilpotent-groups} proves the
nilpotent theorem.

\section{Theorem 1: finite \texorpdfstring{$p$}{p}-groups}
\label{sec:p-groups}

This section compares the vanishing sets of two irreducible characters by
measuring their restriction norms along one fixed chief series. We use
standard notation from character theory; see, for example,
\cite[Chapters~2 and~6]{Isaacs}. If $H\leq G$ and $\alpha,\beta$ are
characters of $H$, then $[\alpha,\beta]_H$ denotes their usual inner
product.

The next lemma isolates, in the notation needed here, the calculation in
\cite[proof of Lemma~3.3]{Wilde}. We recall the argument for completeness;
the calculation itself is due to Wilde.

\begin{lemma}[Wilde's chief-series calculation]\label{lem:chief-jumps}
Let $G\neq 1$ be a finite $p$-group, let
\[
  \mathcal N=(G=N_0\triangleright N_1\triangleright\cdots
  \triangleright N_r=1)
\]
be a chief series, and let $\theta\in\operatorname{Irr}(G)$. Set
\[
  a_i(\theta)=[\theta_{N_i},\theta_{N_i}]_{N_i}
  \qquad (0\leq i\leq r).
\]
For $1\leq i\leq r$, the quotient
\[
  e_i(\theta)=\frac{a_i(\theta)}{a_{i-1}(\theta)}
\]
belongs to $\{1,p\}$, and
\[
  e_i(\theta)=p
  \quad\Longleftrightarrow\quad
  N_{i-1}\setminus N_i\subseteq\operatorname{Van}(\theta).
\]
Consequently,
\begin{equation}\label{eq:chief-product}
  \theta(1)^2
  =\prod_{i=1}^r e_i(\theta)
  =p^{|J_{\mathcal N}(\theta)|},
\end{equation}
where
\[
  J_{\mathcal N}(\theta)
  =\bigl\{i\in\{1,\ldots,r\}:
      N_{i-1}\setminus N_i\subseteq\operatorname{Van}(\theta)\bigr\}.
\]
In particular,
\begin{equation}\label{eq:counting-formula}
  |J_{\mathcal N}(\theta)|=2\log_p\theta(1).
\end{equation}
\end{lemma}

\begin{proof}
We reproduce Wilde's calculation, including the short Clifford-theoretic
point that its restriction norms are powers of $p$. Let $N\trianglelefteq G$
and choose an irreducible constituent $\varphi$ of $\theta_N$. Clifford
theory \cite[Theorem~6.2]{Isaacs} gives
\[
  \theta_N=e(\varphi_1+\cdots+\varphi_t),
  \qquad t=[G:I_G(\varphi)],
  \qquad \theta(1)=et\varphi(1).
\]
The degrees $\theta(1)$ and $\varphi(1)$ are powers of $p$ by
\cite[Theorem~3.11]{Isaacs}, and the index $t$ is a power of $p$.
Consequently $e$ is a power of $p$, while orthogonality of the distinct
constituents $\varphi_1,\ldots,\varphi_t$ gives
\[
  [\theta_N,\theta_N]_N=e^2t.
\]
Thus every $a_i(\theta)$ is a power of $p$.

Fix $i\in\{1,\ldots,r\}$. If $V$ affords $\theta$, then
\[
  [\theta_H,\theta_H]_H
  =\dim_{\mathbb C}\operatorname{End}_H(V)
\]
and
\[
  \operatorname{End}_{N_{i-1}}(V)
  \subseteq\operatorname{End}_{N_i}(V).
\]
Hence $a_{i-1}(\theta)\leq a_i(\theta)$, so
$a_i(\theta)/a_{i-1}(\theta)$ is an integral power of $p$.

Every chief factor of $G$ has order $p$. Therefore
\begin{align}
  p\,a_{i-1}(\theta)
   &=\frac{1}{|N_i|}
     \sum_{x\in N_{i-1}}|\theta(x)|^2 \notag\\
   &=a_i(\theta)
     +\frac{1}{|N_i|}
       \sum_{x\in N_{i-1}\setminus N_i}|\theta(x)|^2.
  \label{eq:norm-comparison}
\end{align}
Every term in the second sum is nonnegative, so
\[
  a_i(\theta)\leq p\,a_{i-1}(\theta).
\]
Combining this inequality with the preceding $p$-power conclusion gives
\[
  e_i(\theta)\in\{1,p\}.
\]

Moreover, equality $a_i(\theta)=p a_{i-1}(\theta)$ holds in
\eqref{eq:norm-comparison} if and only if
\[
  \sum_{x\in N_{i-1}\setminus N_i}|\theta(x)|^2=0.
\]
The summands are nonnegative real numbers. Therefore their sum is zero if
and only if every summand is zero, which is equivalent to
\[
  \theta(x)=0
  \qquad\text{for every }x\in N_{i-1}\setminus N_i.
\]
Thus $e_i(\theta)=p$ exactly when
$N_{i-1}\setminus N_i\subseteq\operatorname{Van}(\theta)$.

Finally,
\[
  a_0(\theta)=[\theta,\theta]_G=1.
\]
The restriction of $\theta$ to the trivial subgroup is $\theta(1)$ times
the trivial character, so
\[
  a_r(\theta)=\theta(1)^2.
\]
Consequently,
\[
  \prod_{i=1}^r\frac{a_i(\theta)}{a_{i-1}(\theta)}
  =\frac{a_r(\theta)}{a_0(\theta)}
  =\theta(1)^2.
\]
The factor $e_i(\theta)$ equals $p$ precisely for
$i\in J_{\mathcal N}(\theta)$ and equals $1$ for every other index. Hence
\[
  \theta(1)^2=p^{|J_{\mathcal N}(\theta)|},
\]
which proves \eqref{eq:chief-product}. Taking exponents of $p$ yields
\eqref{eq:counting-formula}.
\end{proof}

\pgroupdivisibility*

\begin{proof}
The result is immediate when $G=1$, because both irreducible characters
then have degree $1$. Assume that $G\neq1$, and fix a chief series
\[
  \mathcal N=(G=N_0\triangleright N_1\triangleright\cdots
  \triangleright N_r=1).
\]
We use this same series for both characters and compare the indices at which
their restriction norms jump.

Let $i\in J_{\mathcal N}(\chi)$. By the definition of
$J_{\mathcal N}(\chi)$,
\[
  N_{i-1}\setminus N_i
  \subseteq\operatorname{Van}(\chi).
\]
The assumed inclusion of vanishing sets therefore gives
\[
  N_{i-1}\setminus N_i
  \subseteq\operatorname{Van}(\chi)
  \subseteq\operatorname{Van}(\psi).
\]
It follows from the definition of $J_{\mathcal N}(\psi)$ that
$i\in J_{\mathcal N}(\psi)$. Since $i$ was arbitrary,
\[
  J_{\mathcal N}(\chi)\subseteq J_{\mathcal N}(\psi).
\]

By \cite[Theorem~3.11]{Isaacs}, irreducible-character degrees of the
$p$-group $G$ are powers of $p$. Write
\[
  \chi(1)=p^a
  \qquad\text{and}\qquad
  \psi(1)=p^b
\]
for nonnegative integers $a$ and $b$. Applying
\eqref{eq:counting-formula} to the same chief series gives
\[
  2a=|J_{\mathcal N}(\chi)|
  \leq |J_{\mathcal N}(\psi)|=2b.
\]
Hence $a\leq b$, and therefore
\[
  \chi(1)=p^a\mid p^b=\psi(1).
\]

If $\operatorname{Van}(\chi)=\operatorname{Van}(\psi)$, then the inclusion
in each direction gives
\[
  \chi(1)\mid\psi(1)
  \qquad\text{and}\qquad
  \psi(1)\mid\chi(1).
\]
Both degrees are positive integers, so these two divisibilities imply
$\chi(1)=\psi(1)$.
\end{proof}

\section{Theorem 2: finite nilpotent groups}
\label{sec:nilpotent-groups}

We now transfer the $p$-group divisibility theorem to the Sylow direct
factors of a finite nilpotent group.

\nilpotentdivisibility*

\begin{proof}
The result is immediate when $G=1$, because both irreducible characters
then have degree $1$. Assume that $G\neq1$. Since $G$ is finite and
nilpotent, it is the direct product of its Sylow subgroups:
\[
  G=\prod_{p\mid |G|}P_p.
\]
By the direct-product character theorem \cite[Theorem~4.21]{Isaacs},
iterated over these direct factors, there are unique decompositions
\[
  \chi=\bigotimes_{p\mid |G|}\chi_p,
  \qquad
  \psi=\bigotimes_{p\mid |G|}\psi_p,
\]
with $\chi_p,\psi_p\in\operatorname{Irr}(P_p)$.

Fix a prime $p\mid |G|$ and an element
$x\in\operatorname{Van}(\chi_p)$. Regard $x$ as the element
$\widetilde{x}\in G$ whose $p$-coordinate is $x$ and whose other Sylow
coordinates are the identity. The tensor-product decompositions give
\[
  \chi(\widetilde{x})
  =\chi_p(x)\prod_{\substack{q\mid |G|\\q\neq p}}\chi_q(1)
  =0.
\]
Thus $\widetilde{x}\in\operatorname{Van}(\chi)$. By the assumed inclusion,
\[
  \widetilde{x}\in\operatorname{Van}(\psi),
\]
and hence
\[
  0=\psi(\widetilde{x})
  =\psi_p(x)\prod_{\substack{q\mid |G|\\q\neq p}}\psi_q(1).
\]
Every degree $\psi_q(1)$ is a positive integer, so the product over
$q\neq p$ is nonzero. Therefore $\psi_p(x)=0$, which means that
$x\in\operatorname{Van}(\psi_p)$. Since $x$ was arbitrary, we have proved
\[
  \operatorname{Van}(\chi_p)\subseteq\operatorname{Van}(\psi_p).
\]

Theorem~\ref{thm:p-group-divisibility}, applied to $P_p$, now gives
\[
  \chi_p(1)\mid\psi_p(1)
\]
for every prime $p\mid |G|$. Multiplying these divisibilities over the
Sylow direct factors yields
\[
  \chi(1)=\prod_{p\mid |G|}\chi_p(1)
  \mid
  \prod_{p\mid |G|}\psi_p(1)=\psi(1).
\]

If $\operatorname{Van}(\chi)=\operatorname{Van}(\psi)$, applying the
divisibility statement to the two inclusions gives
\[
  \chi(1)\mid\psi(1)
  \qquad\text{and}\qquad
  \psi(1)\mid\chi(1).
\]
The two positive degrees are therefore equal.
\end{proof}

\section*{Declaration of AI use}
OpenAI's ChatGPT was used during the development of this work for assistance
with literature searches, exploration and verification of mathematical arguments,
and manuscript preparation.  All mathematical claims, proofs, citations, and the
final text were independently checked and approved by the author, who takes full
responsibility for the work.

\end{document}